\documentclass{elsart}

\usepackage{natbib}

\usepackage{graphicx}

\usepackage[vcentermath, enableskew]{youngtab}
\newcommand{\Dyoung}[1]{{\tiny \young#1}}
\newcommand{\Dyng}[1]{{\tiny \yng(#1)}}

\usepackage{amsmath}
\usepackage{amssymb}

\usepackage{amsfonts}

\newcommand{\N}{\mathbb{N}}

\newcommand{\Z}{\mathbb{Z}}

\newcommand{\la}{\lambda}
\newcommand{\La}{\Lambda}

\newcommand{\cR}{\mathcal R}
\newcommand{\cQ}{\mathcal Q}
\newcommand{\cH}{\mathcal H}

\newcommand{\0}{{\bf 0}}
\newcommand{\1}{{\bf 1}}

\newcommand{\st}{\;\big | \;} 
\newcommand{\sh}{\mathrm{sh}} 

\begin{document}
\begin{frontmatter}



\title{The Computational Complexity of Rules for the Character Table of $S_n$}


\author{Dan Bernstein}
\ead{dan.bernstein@weizmann.ac.il}
\address{Department of Mathematics,
    The Weizmann Institute of Science,
    Rehovot 76100,
    Israel}

\begin{abstract}
The Murnaghan-Nakayama rule is the classical formula for computing
the character table of $S_{n}$. Y.~Roichman (Roichman 1997) has
recently discovered a rule for the Kazhdan-Lusztig characters of
$q$ Hecke algebras of type $A$, which can also be used for the
character table of $S_{n}$. For each of the two rules, we give an
algorithm for computing entries in the character table of $S_{n}$.
We then analyze the computational complexity of the two
algorithms, and in the case of characters indexed by partitions in
the $(k,\ell$) hook, compare their complexities to each other. It
turns out that the algorithm based on the Murnaghan-Nakayama rule
requires far less operations than the other algorithm. We note the
algorithms' complexities' relation to two enumeration problems of
Young diagrams and Young tableaux.
\end{abstract}

\end{frontmatter}


\DeclareGraphicsExtensions{.eps, .jpg}

\section{Introduction}

This paper examines two formulas for computing entries in the
character table (hereafter called {\it character values}) of the
symmetric group, $S_n$, from the standpoint of computational
complexity. The formulas that we consider are the classical
Murnaghan-Nakayama rule \citep{mur, nak} and the rule recently
discovered by Roichman \citep{roi} for the Kazhdan-Lusztig
characters of $q$ Hecke algebras of type $A$. The discussion is
motivated by a remark in \citet{BR}, in which the authors state
that they are unaware of a comparison of the two rules in terms of
algorithmic complexity, and that ``one would expect that they have
the same complexity".

The irreducible characters of $S_{n}$ are a distinguished set of
class functions $\{\chi^{\la}:S_n \to \Z \st \la \vdash n\}$
\citep[for a complete description see][]{sag}. A character value
of $S_n$ is indexed by an ordered pair $(\la,\mu)$ of partitions
of $n$ and is denoted by
\[\chi^\la(\mu) = \chi^\la(w) \qquad w\in S_n \text{ is of cycle
type $\mu$}.\]

A formula for $\chi^{\la}(\mu)$ suggests a systematic way for
computing character values --- an algorithm whose input is a pair
of partitions $(\la,\mu)$ and whose output is the integer
$\chi^\la(\mu)$. It is such an algorithm's computational
complexity that is examined for each of the two rules.

The rest of the paper is organized as follows: in section
\ref{SEC:MN} we present the Murnaghan-Nakayama rule and specify an
algorithm based on it. In section \ref{SEC:Roi} we treat
Roichman's rule similarly.

Section \ref{SEC:instance} gives the complexity of computing a
single character value using each of the two algorithms. Two
enumeration problems, of Young diagrams and of Young tableaux,
occur in the dominant factors in the complexity of the
Murnaghan-Nakayama rule (eq.~\ref{EQ:tMN}) and Roichman's rule
(eq.~\ref{EQ:tRoi}), respectively.

Finally, in section \ref{SEC:comp}, we compare the algorithms in
terms of their worst-case complexity on the family of characters
indexed by partitions in the $(k,\ell)$ hook.

\subsection{Main Results}

Given $n$ and partitions $\la$ and $\mu$ of $n$ we show that:
\begin{enumerate}
\item
    The running time of our Murnaghan-Nakayama-based algorithm is,
    up to a factor of order $n$, the number of Young diagrams
    that are contained in the Young diagram of $\la$ and satisfy
    an additional constraint determined by $\mu$ (see proposition
    \ref{PR:tMN1}). Lemma \ref{LE:r_la} gives a determinantal
    formula for this number when the constraint is empty.

\item
    The running time of the algorithm based on Roichman's rule is,
    up to a factor of order $n$, the number of standard Young
    tableaux whose shape is contained in the Young diagram of
    $\la$ that satisfy an additional constraint determined by
    $\mu$ (see proposition \ref{PR:tRoi}). By lemma \ref{LE:q_la},
    when the constraint is empty, this number is $O(n)d_\la$,
    where $d_\la$ is the degree of $\la$, that is the number of standard Young tableaux of
    shape $\la$.
\end{enumerate}

Given $k$ and $\ell$, worst-case analysis of the family of
characters where the choice of $\la$ is restricted to the
$(k,\ell)$ hook shows that in this case the
Murnaghan-Nakayama-based algorithm's complexity is
$\Theta(n^{k+\ell+1})$ (see proposition \ref{PR:klhookMN}) whereas
the complexity of the algorithm based on Roichman's rule is in
$\Omega(n^{-g} (k+\ell)^n) \cap O(n^{-g+2} (k+\ell)^n)$ for some
constant $g$ (see proposition \ref{PR:klhookRoi}).

Some experimental results for characters not in the above family
are also included in subsection \ref{SUB:GeneralDiagrams}.

\section{The Murnaghan-Nakayama Rule}\label{SEC:MN}

A shape $A\subset \N\times\N$ is said to be {\it edgewise
connected} if \[A=\{(i_{1},j_{1}),
(i_{2},j_{2}),\dots,(i_{n},j_{n})\}\] and for all $k<n$ \[| i_{k}-
i_{k+1} | + |j_{k} - j_{k+1}| = 1\] (i.e. each cell is exactly one
horizontal or vertical step away from its predecessor). For
example, $\Dyoung{(\hfil\hfil,\hfil)}$ is edgewise connected, but
$\Dyoung{(:\hfil,\hfil)}$ is not.

Let $\la=(\la_1,\la_2,\dots,\la_m) \vdash n$. A skew diagram
$\xi=\la / \mu$ is said to be a {\it rim hook} of $\la$ if $\xi$
is edgewise connected and contains no $2\times 2$ subset of cells
(\Dyng{2,2}). In this case we write $\la \backslash \xi = \mu$ and
say that $\mu$ is obtained by {\it removing} the rim hook $\xi$
from $\la$. For example, if $\la=(4,3,2)$, then \[\la / (2,2,2) =
\young(::\hfil\hfil,::\hfil)\] is a rim hook of $\la$, but \[\la /
(2,2,1) = \young(::\hfil\hfil,::\hfil,:\hfil) \text{ and } \la /
(1,1) = \young(:\hfil\hfil\hfil,:\hfil\hfil,\hfil\hfil)\] are not:
the former is not edgewise connected, and the latter contains a
$2\times 2$ block.

The {\it leg length of a rim hook} $\xi$ is \[ll(\xi) = (\text{the
number of rows of $\xi$})-1.\]

Let $\la \backslash \la_{1}$ denote the partition
$(\la_{2},\la_{3},\dots,\la_{m})$.

Note that the notation $\la / \mu$ is reserved for skew diagrams,
while $\la \backslash \xi$ and $\la \backslash \la_1$ are always
ordinary diagrams.
\newpage
The following is the classical recursive formula for computing
characters of $S_{n}$.

\begin{thm}[The Murnaghan-Nakayama Rule] \label{T:MN}
Let $\la, \mu \vdash n$. Then
\begin{equation}\label{EQ:MN}
\chi^\la(\mu) = \sum_{\xi}(-1)^{ll(\xi)} \chi^{\la\backslash\xi}(\mu\backslash\mu_{1})
\end{equation}
where the sums runs over all rim hooks $\xi$ of $\la$ having
$\mu_{1}$ cells, and $\chi^{(0)}(0)=1$.
\end{thm}
A proof appears in \citet{sag}.

\begin{exmp}\label{EX:MN}
Calculating $\chi^\la(\mu)$ where $\la=(5,4,2,1)$ and
$\mu=(4,3,2,2,1)$.

The computation process can be viewed as a tree. The appropriate
signs appear beside the arrows indicating the removal of rim
hooks.

\[
    \begin{array}{l|cccccccc}
({\bf4},3,2,2,1)&   & & & & & \Dyng{5,4,2,1} & & \\
&   & & & & -\swarrow & & \searrow- & \\
({\bf3},2,2,1)& & & \Dyng{5,1,1,1} & & & & & \Dyng{3,2,2,1} \\
&   & +\swarrow & & \searrow+ & & & & -\downarrow\quad \\
({\bf2},2,1)&   \Dyng{5} & & & & \Dyng{2,1,1,1} & & & \Dyng{3,2} \\
&   +\downarrow\quad & & & &-\downarrow\quad & & & +\downarrow\quad \\
({\bf2},1)& \Dyng{3} & & & & \Dyng{2,1} & & & \Dyng{3} \\
&   +\downarrow\quad    & & & & 0           & & & +\downarrow\quad \\
({\bf1})&   \Dyng{1} & & & &                & & & \Dyng{1} \\
&   +\downarrow \quad   & & & &                 & & & +\downarrow\quad \\
()& \emptyset   & & & &             & & & \emptyset \\
&       1       & & & &             & & &   1
    \end{array}
\]
\begin{align*}
\chi^{(5,4,2,1)}(4,3,2,2,1) & = -\chi^{(5,1,1,1)}(3,2,2,1)-\chi^{(3,2,2,1)}(3,2,2,1) \\
& = -(\chi^{(5)}(2,2,1)+\chi^{(2,1,1,1)}(2,2,1)) +\chi^{(3,2)}(2,2,1) \\
& = -(\chi^{(3)}(2,1)-\chi^{(2,1)}(2,1))+\chi^{(3)}(2,1) \\
& = -(\chi^{(1)}(1)-0)+\chi^{(1)}(1) \\
& = -(1+0)+1 \\ &= 0.
\end{align*}

\end{exmp}

\subsection{An algorithm based on the Murnaghan-Nakayama rule}

Computing the sum in \eqref{EQ:MN} requires enumerating all rim
hooks of certain length of a given partition. This is done using
partition sequences (\citet{ols}, \citet{bes}).

A {\it partition sequence} $\La$ is a doubly infinite sequence of
binary digits starting with an infinite sequence of zeros and
ending with an infinite sequence of ones. For example, \[\La=
\dots \;\0\;\0\;\underbrace{ \1\;\0\;\1\;\0\;\1\;\1\;\0\;\1\;\0
}_{\bar\La} \;\1\;\1\; \dots \] where the dots at the beginning
(end) represent an infinite sequence of \0s (\1s), is a partition
sequence.

We shall refer to the finite subsequence of $\La$ starting with
the first \1 and ending with the last \0 as the {\it essential
part of $\La$}, which we will denote by $\bar \La$.

Given a partition $\la=(\la_1,\la_2,\dots,\la_m)$, its partition
sequence is defined as
\[
\La = \dots \0 \;\0\; \underbrace{\1\; \1 \;\dots \;\1}_{\la_{m}} \;\0 \; \underbrace{\1\; \1 \;\dots \;\1}_{\la_{m-1}-\la_{m}} \0\;\underbrace{\1\; \1 \;\dots \;\1}_{\la_{m-2}-\la_{m-1}} \0\; \1\; \dots \;\1 \;\0\; \underbrace{\1\; \1 \;\dots \;\1}_{\la_{2}-\la_{1}}\; \0 \;\1\;\1\; \dots
\]

For example, the partition sequence of $\la=(5,4,2,1)$ is \[\La =
\dots \;\0\;\0\;  \1\;\0\;\1\;\0\;\1\;\1\;\0\;\1\;\0   \;\1\;\1\;
\ldots\]

The graphic version of this construction is a walk along the
borderline of $\la$, coming form the south on the vertical line,
going along the border and leaving on the horizontal line
eastwards, recording each vertical step by a \0 and each
horizontal step by a \1. In our example, the borderline of the
Young diagram of $\la=(5,4,2,1)$ is \[ \includegraphics{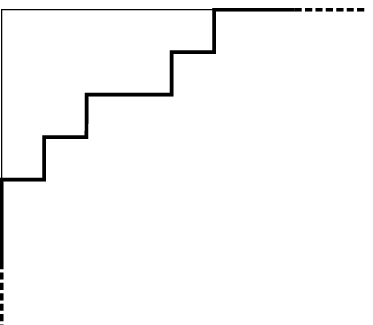}
\] which indeed gives the sequence
\[\dots \;\0\;\0\;  \1\;\0\;\1\;\0\;\1\;\1\;\0\;\1\;\0
\;\1\;\1\; \ldots = \La .\]

Consider the rim hook $\xi$ of $\la=(5,4,2,1)$:
\[
\includegraphics{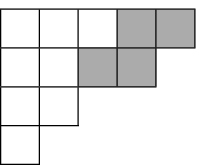} 
.\]
The partition sequence of $\mu = \la \backslash \xi =
(3,2,2,1)$ is
\[
\begin{array}{rcccccccccccccc}
M=\dots & \0 & \0 & \1 & \0 & \1 & \0 & {\fbox \0} & \1 & \0 & \1
& {\fbox \1} & \1 & \1 & \dots
\end{array}
\]
and the partition sequence of $\la$ is
\[
\begin{array}{rcccccccccccccc}
\La=\dots & \0 & \0 & \1 & \0 & \1 & \0 & {\fbox \1} & \1 & \0 & \1 & {\fbox \0} & \1 & \1 & \dots \\
\end{array}
\]
We observe that $M$ differs from $\La$ only by the exchange with one another of the two digits in the positions marked in the above sequences, changing their order from the \1 being to the left of the \0 in $\La$, to the \0 being to the left of the \1 in $M$. Moreover, we note that the two digits are $4=|\xi|$ positions apart from each other and that there is exactly $1=ll(\xi)$ \0 between them.

This is not by coincidence, as the following definitions and the
next proposition show.

Let $\La$ be a partition sequence. A {\it rim hook in} $\La$ is a
pair consisting of a \0 and a \1 in $\La$ such that the \1 appears
to the left of the \0.
The distance between the \0 and the \1 is
the {\it length} of the rim hook, and the number of \0s strictly
between them is the {\it leg length} of the rim hook. The rim hook
is {\it removed} by exchanging the \0 with the \1. For example,
the marked pair of digits in
\[
\begin{array}{rcccccccccccccc}
\La=\dots & \0 & \0 & \1 & \0 & \fbox{\1} & \0 & \1 & \1 & \0 & \1 & \fbox{\0} & \1 & \1 & \dots \\
\end{array}
\]
is a rim hook of length 6 and of leg length 2. The partition sequence obtained from $\La$ by removing this rim hook is
\[
\begin{array}{rcccccccccccccc}
\qquad\dots & \0 & \0 & \1 & \0 & \fbox{\0} & \0 & \1 & \1 & \0 & \1 & \fbox{\1} & \1 & \1 & \dots \\
\end{array}
\]

\begin{prop}
There is a bijection between rim hooks in the partition sequence
$\La$ of a partition $\la$ and rim hooks of the Young diagram of
$\la$. Moreover, this bijections preserves the notions of length
and leg length, and the removal of a rim hook in the partition
sequence corresponds to the removal of the corresponding rim hook
of the Young diagram.
\end{prop}

Based on this, the following algorithm, {\it MNinner}, finds all
rim hooks $\xi$ of $\la$ having length $\mu_{1}$ simply by going
over all pairs of digits that are $\mu_{1}$ places apart from each
other in $\bar\La$, where $\La$ is the partition sequence of
$\la$. A variable $\sigma$ keeps track of $(-1)^{\text{\# of \0s
between the 2 digits of the pair}}$. If and only if the left digit
in such a pair is \1 and the right digit is \0, then it is a rim
hook, and then the partition sequence of $\la\backslash\xi$ is
obtained by exchanging the \0 with the \1. {\it MNinner} then
proceeds recursively to compute
$\chi^{\la\backslash\xi}(\mu\backslash\mu_{1})$ and adds
$\sigma\chi^{\la\backslash\xi}(\mu\backslash\mu_{1})$ to the sum.

\begin{tabbing} {\sf Fun}\={\sf ction}  {\it MNinner} ($R, t$) \\
{\sf Input:} A sequence of binary digits $R=R_{1}R_{2}\dots R_{s}$ and an index $t$ \\
{\sf Out}\={\sf put:} \=$\chi^{\rho}(\nu)$ where $\rho$ is the partition whose partition sequence\\
        \>\>is $\dots\;\0\;R\;\1\;\dots$ and $\nu=(\mu_{t},\mu_{t+1},\dots,\mu_{k})$. \+ \\
    {\sf if} \= $t>k$ \+\\
    {\sf then} $\chi\leftarrow 1$
    \\
    {\sf else} \= \+ $\chi \leftarrow 0$ \`
    \\
    $\sigma \leftarrow 1$ \\
     {\sf for} \= $j \leftarrow 1$ {\sf to} $\mu_t-1$ \+ \`
     \\
    {\sf do} {\sf if} $R_j=0$ {\sf then} $\sigma \leftarrow -\sigma$ \- \\
    {\sf for} \= $i \leftarrow 1$ {\sf to} $s-\mu_t$\+ \\
    {\sf do} \= \+ {\sf if} $R_{i} \neq R_{i+\mu_t-1}$ {\sf then} $\sigma \leftarrow -\sigma$
    \\
        {\sf if} \=  $(R_i,R_{i+\mu_t})$ is a rim hook, $\xi$ \+ \\
        {\sf then} \=\+ exchange between $R_{i}$ and $R_{i+\mu_{t}}$. \\
          $\chi \leftarrow \chi+\sigma \cdot${\it MNinner}($R,t+1$)\\
          exchange between $R_{i}$ and $R_{i+\mu_{t}}$. \- \-  \- \- \- \-\\
{\sf return} $\chi$
\end{tabbing}

A major inefficiency of {\it MNinner} is that character values
that occur more than once in the expansion of the right hand side
of \eqref{EQ:MN} are re-computed each time. In example
\ref{EX:MN}, $\chi^{(3)}(2,1)$ occurs twice (and therefore so does
$\chi^{(1)}(1)$), so {\it MNinner} is invoked twice to compute it.
This is overcome in the following algorithm, {\it MN1inner}, by
saving intermediate results in a table and using it to look up
character values before computing them. Each time a value
$\chi^{\rho}(\nu)$ is computed, it is recorded in a table $T$, and
each time a value is required, it is first looked up in $T$, and
only if it is not there, then it is computed. $T$ is indexed by
partitions $\rho$ rather than by pairs $(\rho, \nu)$ of
partitions, since for any $\chi^{\rho}(\nu)$ appearing in the
expansion of \eqref{EQ:MN} we have that $\nu$ is the tail of $\mu$
of weight $|\rho|$.
\newpage
\begin{tabbing} {\sf Fun}\={\sf ction}  {\it MN1inner} ($R, t$) \\
{\sf Input:} A sequence of binary digits $R=R_{1}R_{2}\dots R_{s}$ and an index $t$ \\
{\sf Out}\={\sf put:} \=$\chi^{\rho}(\nu)$ where $\rho$ is the partition whose partition sequence\\
        \>\>is $\dots\;\0\;R\;\1\;\dots$ and $\nu=(\mu_{t},\mu_{t+1},\dots,\mu_{k})$. \+ \\
    {\sf if} \= $t>k$ \+\\
    {\sf then} $\chi\leftarrow 1$
    \\
    {\sf else} \= \+ $\chi \leftarrow 0$ \`
    \\
    $\sigma \leftarrow 1$ \\
     {\sf for} \= $j \leftarrow 1$ {\sf to} $\mu_t-1$ \+ \`
     \\
    {\sf do} {\sf if} $R_j=0$ {\sf then} $\sigma \leftarrow -\sigma$ \- \\
    {\sf for} \= $i \leftarrow 1$ {\sf to} $s-\mu_t$\+ \\
    {\sf do} \= \+ {\sf if} $R_{i} \neq R_{i+\mu_t-1}$ {\sf then} $\sigma \leftarrow -\sigma$
    \\
        {\sf if} \=  $(R_i,R_{i+\mu_t})$ is a rim hook, $\xi$ \+ \\
        {\sf then} \=\+ exchange between $R_{i}$ and $R_{i+\mu_{t}}$. \\
        {\sf if} \= \+ $T(\tilde\rho)$ is empty, where $\tilde\rho=\rho\backslash\xi$ is the partition\\
        whose partition sequence is \dots\0\;$R$\;\1\dots\\
          {\sf then} \=\+ $T(\tilde\rho)\leftarrow${\it MN1inner}($R,t+1$) \-\-\\

          $\chi \leftarrow \chi+\sigma T(\tilde\rho)$\\
          exchange between $R_{i}$ and $R_{i+\mu_{t}}$. \- \- \- \- \- \-\\
{\sf return} $\chi$
\end{tabbing}

Given partitions $\la$ and $\mu$, to compute $\chi^{\la}(\mu)$ one
needs to compute the essential part of $\la$'s partition sequence,
$\bar\La$, and then to invoke {\it MN1inner}($\bar\La,1$). This is
what algorithms {\it PartSeq} and {\it MurNak} do.

\begin{tabbing}
    {\sf Fun}\={\sf ction} {\it PartSeq}($\la$) \\
    {\sf Input:} a partition $\la=(\la_{1},\la_{2},\dots,\la_{m})$. $m=\ell(\la)$ \\
    {\sf Output:} $\bar\La$, the essential part of the partition sequence of $\la$ \+\\
        $\bar\La \leftarrow$ an empty sequence\\
        $\la_{m+1} \leftarrow 0$\\
        {\sf for} \= $i \leftarrow m$ {\sf down to} $1$\+  \\
            {\sf do} \=\+ {\sf for} \=\+ $k \leftarrow 1$ {\sf to}  $\la_{i}-\la_{i+1}$\\
             {\sf do} \=\+ $\bar\La \leftarrow \bar\La \| 1$ \-\- \\
            $\bar\La \leftarrow \bar\La \| 0$\- \- \\
   {\sf return $\bar\La$}
\end{tabbing}

\begin{tabbing}
{\sf Fun}\={\sf ction}  {\it MurNak} ($\la,\mu$) \\
{\sf Input:} partitions $\la$ and $\mu$ of the same weight\\
{\sf Output:} $\chi^{\la}(\mu)$\+ \\
    $T \leftarrow$ a 1-dimensional sparse array \\
    $\bar\La \leftarrow \textit{PartSeq}(\la)$\\
    $\chi \leftarrow \textit{MN1inner}(\bar\La,1)$ \\
{\sf return} $\chi$
\end{tabbing}

\section{Roichman's Rule}\label{SEC:Roi}

Let $(W,S)$ be a Coxeter system, and let $\ell(w)$, $w \in W$ be
the length function with respect to $S$.

The $q$ Hecke algebra $\cH$ of $W$ is the algebra spanned by the
set $\{ T_w \st w \in W \}$ over the ring of Laurent polynomials
$\Z [q, q^{-1}]$ subject only to the following relations:
\[
T_s T_w = T_{sw} \quad \text{if $s \in S$ and $\ell(sw)>\ell(w)$ }
\]
\[
T_s^2 = (q-1)T_s + qT_1 \quad \text{if $s \in S$}
\]
where $T_1$ acts as the identity.

\citet{KL} gives a distinguished basis $\{C_w \st w\in W \}$ for
$\cH$ and a partition of the Coxeter group $W$ into {\it
Kazhdan-Lusztig cells}. Each left Kazhdan-Lusztig cell $\mathcal
C$ has a left representation of $\cH$ associated to it. Let
$\chi^{\mathcal C}$ be the character of that representation. Then
for any $T\in\cH$ and finite cell $\mathcal C$
\begin{equation}\label{EQ:chi^C(T)}
\chi^{\mathcal C}(T) = \sum_{w\in \mathcal C} TC_w(w)
\end{equation}
where $TC_w(w)$ is the coefficient of $C_w$ in $TC_w$.

\citet{roi} gives a formula for $T_{s_1 s_2 \dots s_k}C_w(w)$
where $s_1,\dots,s_k \in S$, subject to certain relations between
the $s_i$. In the case $W=S_n$, the formula applies to all of the
summands in \eqref{EQ:chi^C(T)}. Furthermore in the $W=S_n$ case,
the Kazhdan-Lusztig characters are exactly the irreducible
characters, and the Robinson-Schensted-Knuth correspondence gives
rise to a canonical map between the Kazhdan-Lusztig and Young's
natural characters of $S_n$, allowing for the formulation of the
character as a weighted sum over standard tableaux.

If $\mu=(\mu_1,\dots,\mu_k) \vdash n$, define
$B(\mu)=\{\mu_{1}+\dots+\mu_{r} \st 1\le r \le k\}$. For example,
if $\mu=(5,2,1,1)$ then $B(\mu)=\{5, 7, 8, 9\}$.

Recall that a {\it standard tableau} is a tableau whose rows and
columns are increasing. The {\it descent set} of a standard
tableau $T$ is
\[
D(T) = D(w(T)^{-1})
\]
the descent set of the inverse of the reading word of $T$, also
characterized by
\begin{align*}
D(T) &= \{i \st i+1 \text{ is in the southwest of } i \text{ in
}T\}
\end{align*} where ``southwest" means strictly below and weakly to the left.
For example, the descent set of $T=\Dyoung{(12458,36,79)}$ is
$\{2,5,6,8\}$.

Define $f^q_{\mu}(T,i)$, $i=1,2,\dots,n$ by
\begin{equation}\label{EQ:fmuQi}
f^q_{\mu}(T,i)=\begin{cases}
    -1  &\text{$i\notin B(\mu)$, $i\in D(T)$} \\
    0   &\text{$i,i+1\notin B(\mu)$, $i\notin D(T)$ and $i+1\in D(T)$} \\
    q   &\text{otherwise} .
\end{cases}
\end{equation}

\begin{thm}[Roichman's Rule \citep{roi}]\label{T:Roi}
Let $\la,\mu\vdash n$, and let $\chi^\la$ be the corresponding
character of the $q$ Hecke algebra of $S_n$. Let $T_\mu$ be any
element in the Hecke algebra indexed by a permutation $w\in S_n$
of cycle type $\mu$. Then
\[
\chi^\la(T_\mu)=\sum_T \prod_{\substack{1\le i <n \\
i\notin B(\mu)}}f^q_{\mu}(T,i)
\]
where the sum runs over all standard tableaux $T$ of shape $\la$.
\end{thm}

Substituting 1 for $q$ in the above we get a rule for the
characters of $S_n$:
\begin{equation}\label{EQ:Roi}
\chi^\la(\mu)=\sum_T \prod_{1\le i <n}f^1_{\mu}(T,i)
\end{equation}
where the sum runs over all standard tableaux $T$ of shape $\la$.

\begin{exmp}\label{EX:Roi}
Calculating $\chi^\la(\mu)$ where $\la=(2,1,1)$ and $\mu=(3,1)$.
We have $B(\mu)=\{3,4\}$ and
\[
\begin{tabular}{|c|c|c|c|c|}
\hline \parbox[c][20pt]{1pt}{$T$} & $f^1_\mu(T,1)$ & $f^1_\mu(T,2)$ & $f^1_\mu(T,3)$ & $\prod_{1\le i<4}f^1_\mu(T,i)$ \\
\hline \parbox[c][36pt]{15pt}{\Dyoung{(12,3,4)}} &  $0$ & $-1$ & $1$ &  $0$ \\
\hline \parbox[c][36pt]{15pt}{\Dyoung{(13,2,4)}} & $-1$ &  $1$ & $1$ & $-1$ \\
\hline \parbox[c][36pt]{15pt}{\Dyoung{(14,2,3)}} & $-1$ & $-1$ &
$1$ &  $1$ \\
\hline
\end{tabular}
\]
Hence $\chi^\la(\mu) = 0+(-1)+1 = 0$.
\end{exmp}

\subsection{Recursive formulation}

The na\"{\i}ve way to compute $\chi^\la(\mu)$ using Roichman's
rule would be to construct all standard tableaux of shape $\la$,
and for each tableau $T$ to compute the values of $f^1_{\mu}(T,i)$
for all $i$ (or until a 0 value is encountered), and finally to
take the products and their sum. However, it can easily be shown
that $f^q_{\mu}(T,i)$ depends only on the first $i+2$ entries of
$T$. An improvement over the na\"{\i}ve approach is achieved by
using this observation, as follows:

Let $T$ be a standard tableau. Denote by $T_{|i}$ the standard
tableau obtained by deleting all entries $j>i$ from $T$. For
example, $\Dyoung{(1246,357)}_{|4}=\Dyoung{(124,3)}$.

Define
\[
g_\mu(T,i) = \begin{cases}
        1   &   i=1,2   \\
        f^1_\mu(T,i-2)    &   i>2 .
    \end{cases}
\]

Let $Q$ be a standard tableau of shape $\alpha\subseteq\la$,
$|\alpha|=j$. Define
\begin{equation}\label{EQ:A}
A(\la,\mu,Q)= \sum_{T}\prod_{j<i<n+2}g_{\mu}(T_{|i},i)
\end{equation}
where the sum runs over all standard tableaux $T$ of shape $\la$
containing $Q$ (i.e. such that $T_{|j}=Q$). Note that
\begin{align*}
A(\la,\mu,\emptyset) &=  \sum_{T}\prod_{0<i<n+2}g_{\mu}(T_{|i},i)\\
    &= \sum_{T}\prod_{0<i<n+2}g_{\mu}(T_{|i},i)\\
    &= \sum_{T}\prod_{2<i<n+2}f^1_{\mu}(T_{|i},i-2)\\
    &= \sum_{T}\prod_{0\le i<n}f^1_{\mu}(T_{|i+2},i)\\
    &= \sum_{T}\prod_{0\le i<n}f^1_{\mu}(T,i)
\end{align*}
where the sums run over all standard tableaux of shape $T$, so
\eqref{EQ:Roi} can be rewritten as:
\begin{equation}\label{EQ:chiEqualsA}
    \chi^\la(\mu)=A(\la,\mu,\emptyset).
\end{equation}

The following proposition follows easily from the definitions.
\begin{prop}\label{PR:recursiveA}
We have the following recursive formula:
\begin{equation}\label{EQ:recursiveA}
A(\la,\mu,Q)= \begin{cases}
    g_{\mu}(Q,n+1)  &j=n \\
    \sum_{S}g_{\mu}(S,j+1)A(\la,\mu,S)  &\text{otherwise}
    \end{cases}
\end{equation}
where the sum runs over all standard tableaux $S$ such that
$\sh(S)\subseteq \la$, $Q\subset S$ and $|S|=|Q|+1$, i.e. the
tableaux $S$ are those obtained by adding $j+1$ to $Q$ in a
position belonging to $\la$.
\end{prop}

\begin{exmp}
Calculating $\chi^\la(\mu)$ where $\la=(2,1,1)$ and $\mu=(3,1)$,
hence $B(\mu)=\{3,4\}$. Each node in the following tree shows a
tableau $Q$ of some shape contained in $\la$, starting with the
empty tableau. Each node's children are all the tableaux $S$ of
shape contained in $\la$ that can be obtained by adding one entry
to $Q$. The number in parentheses above each tableau $Q$ is
$g_{\mu}(Q,j)$ where $j$ is the number of entries in $Q$, also
appearing in the column to the left of the tree. The numbers in
the last row are $g_\mu(\cdot,5)$ for the tableaux above them.
\[
\begin{array}{l|ccccccc}
&   &   &\Dyoung{(\hfil\hfil,\hfil,\hfil)} \\
&   &   &\downarrow \\
1&  &   &\substack{(1) \\ \Dyoung{(1\hfil,\hfil,\hfil)}} \\
&   &\swarrow   &       & \searrow \\
2&\substack{(1)\\ \Dyoung{(12,\hfil,\hfil)}} & &   &   & \substack{(1)\\ \Dyoung{(1\hfil,2,\hfil)}} \\
&\downarrow &   &   & \swarrow  &   & \searrow \\
3&\substack{(0)\\ \Dyoung{(12,3,\hfil)}} & & & \substack{(-1)\\ \Dyoung{(13,2,\hfil)}} & & \substack{(-1)\\ \Dyoung{(1\hfil,2,3)}} \\
& \downarrow  & & &  \downarrow & & \downarrow \\
4& \substack{(-1)\\ \Dyoung{(12,3,4)}} & & & \substack{(1)\\ \Dyoung{(13,2,4)}} & & \substack{(-1)\\ \Dyoung{(14,2,3)}} \\
\\
5 & 1 &  & & 1 & & 1
\end{array}
\]
$\chi^{(2,1,1)}(3,1)=1\left(1\cdot 0 + 1\left( -1\cdot 1 \cdot 1 +
(-1)(-1)1\right)\right) = 0$.
\end{exmp}

\subsection{An algorithm based on Roichman's rule}

The following algorithm, {\it RoiInner}, computes $A(\la,\mu,Q)$
according to proposition \ref{PR:recursiveA}, computing values of
$g_\mu$ and invoking itself recursively as necessary. It assumes
that the global variable $B=B_{1}B_{2}\dots B_{n}$ is assigned the
values $B_{m}=1_{B(\mu)}(m)$. In the case $j=n$ of
\eqref{EQ:recursiveA}, it computes
$g_{\mu}(Q,n+1)=f^1_{\mu}(Q,n-1)_{|q=1}$ by checking for the first
case of \eqref{EQ:fmuQi} (note that the second case of
\eqref{EQ:fmuQi} cannot occur for $i=n-1$ since $n\in B(\mu)$
always). In the case $j<n$ of \eqref{EQ:recursiveA}, the algorithm
computes the sum by going over the rows of $Q$, checking for each
row whether by adding $j+1$ at its end one gets a tableau $S$ such
that $\sh(S)\subset \la$. If so, it sets the variable $\tilde d$
to indicate whether $j\in D(S)$ and determines $g_{\mu}(S,j+1)$,
which is assigned to the variable $g$. Finally, if $g\neq 0$, it
proceeds recursively to compute $A(\la,\mu,S)$ and adds
$g_{\mu}(S,j+1)A(\la,\mu,S)$ to the sum.

\newpage
\begin{tabbing}
{\sf Fun}\={\sf ction}  {\it RoiInner}($\alpha,j,m,d$) \\
{\sf Input:} a partition $\alpha=(\alpha_{1},\dots,\alpha_{\ell})$, $j=|\alpha|$, a row index $m$ and an indicator $d$. \\
{\sf Out}\={\sf put:} \=$A(\la,\mu,Q)$ where $Q$ is any of the tableaux such that\\
 \>\> $\sh(Q)=\alpha$, $j$ appears on row $m$ of $Q$ and $d$ indicates whether $j-1\in D(Q)$.\+\\
     {\sf if} \= $j=n$\+\\
      {\sf then} \=\+
        {\sf if} \=\+ ($d=\textsf{yes}$ {\sf and} $B_{n-1}=0$)\\
         {\sf  then}\=  $A \leftarrow -1$  \\
        {\sf el}{\sf se} \= $A \leftarrow 1$ \- \-\\
    {\sf el}{\sf se} \=\+ $A\leftarrow 0$ \\
      {\sf for} \=\+ $k \leftarrow 1$ {\sf to} $\ell$\\
     {\sf do} {\sf if} \=\+ (($k=1$ {\sf or} $\alpha_{k}<\alpha_{k-1}$) {\sf and} $\alpha_{k}<\la_{k}$)\\
      {\sf then} \=\+  {\sf if} \=\+ $k>m$ \\
         {\sf then} $\tilde d\leftarrow$ {\sf yes}\\
         {\sf else} $\tilde d \leftarrow$ {\sf no}\- \\
      {\sf if} \=\+ ($j+1<3$ {\sf or} $B_{j-1}=1$)\\
        {\sf then}   $g \leftarrow 1$  \\
        {\sf el}{\sf se if} \=\+ $d=\sf{yes}$ \\
              {\sf then}       $g \leftarrow -1$ \\
              {\sf el}{\sf se if} \=\+ ($\tilde d=${\sf yes and} $d=${\sf no and} $B_j=0$)\\
                {\sf then}  $g \leftarrow 0$  \\
                {\sf el}{\sf se} $g \leftarrow 1$ \-\-\- \\
                {\sf if} \=\+ $g\neq 0$\\
                  {\sf then}  \=\+ $\alpha_k \leftarrow \alpha_k+1$  \\
                    $A \leftarrow A+g\cdot\textit{RoiInner}(\alpha,j+1,k,\tilde d)$ \\
                    $\alpha_k \leftarrow \alpha_k-1$ \- \- \- \- \- \-\- \ \\
{\sf return} A
\end{tabbing}

Given partitions $\la$ and $\mu$, to compute $\chi^{\la}(\mu)$ one
needs to initialize the global variable $B$ to contain $B(\mu)$
and then to compute $A(\la,\mu,\emptyset)$ using {\it RoiInner}.
This is what the algorithm {\it Roich} does.

\begin{tabbing}
{\sf Fun}\={\sf ction}  {\it Roich}($\la,\mu$)  \\
{\sf Input:} partitions $\la=(\la_1,\dots,\la_\ell)$ and $\mu=(\mu_1,\dots,\mu_k)$ of the same weight. $\ell=\ell(\la)$ \\
{\sf Output:} $\chi^{\la}(\mu)$ \+\\
    $\alpha \leftarrow ( \underbrace{0,0,\dots,0}_\ell )$\\
    {\sf for} \= $i \leftarrow 1$ {\sf to} $k$\+\\
     {\sf do} \=\+ {\sf for} \=\+ $j \leftarrow 1$ {\sf to} $\mu_i - 1$\\
      {\sf do} $B \leftarrow B \| 0$ \- \\
       $B \leftarrow B \| 1$ \-\- \\
  $\chi \leftarrow ${\it RoiInner}($\alpha, 0, 1,\sf{no}$) \\
 {\sf return} $\chi$
\end{tabbing}

\section{Problem Instance Complexity}\label{SEC:instance}

A {\it problem instance} in the case of computing character values
of the symmetric group is simply an ordered pair $(\la,\mu)$ of
partitions of the same integer.

\subsection{\it MurNak}

Let $\cR_{\la,\mu}$ be the set of partitions appearing in the
expansion of the right hand side of \eqref{EQ:MN}. More precisely,
define
\begin{multline*}
\cR_{\la,\mu}=\{ \alpha \subseteq \la \st \exists
\alpha=\alpha_{i}\subset \alpha_{i-1}\subset \dots \subset
\alpha_{0}=\la \text{, }\\
    \xi_{j}=\alpha_{j-1}/\alpha_{j} \text{
is a rim hook}, |\xi_{j}|=\mu_{j}\text{, }1\le j\le i\}
\end{multline*}
which is the set of partitions one can obtain starting with $\la$
by removing a sequence of rim hooks $\xi_1,\dots,\xi_i$ of lengths
$\mu_1,\dots,\mu_i$ respectively, $i\le \ell(\mu)$. In example
\ref{EX:MN}, for instance,
\begin{equation*}
\cR_{(5,4,2,1),(4,3,2,2,1)} = \left\{ \Dyng{5,4,2,1},
\Dyng{5,1,1,1}, \Dyng{3,2,2,1}, \Dyng{2,1,1,1}, \Dyng{5},
\Dyng{3,2}, \Dyng{2,1}, \Dyng{3}, \Dyng{1}, \emptyset \right\}.
\end{equation*}

Denote $|\cR_{\la,\mu}|$ by $r_{\la,\mu}$.

\begin{prop}\label{PR:tMN1}
Let $t_{\textit MurNak}(\la,\mu)$ be the running time of {\it
MurNak} on input $(\la,\mu)$. Then
\begin{equation}\label{EQ:tMN}
    t_{\textit MurNak}(\la,\mu) \in \Theta(r_{\la,\mu} h_{1,1}(\la))
\end{equation}
Where $h_{1,1}(\la) = \la_{1}+\la'_{1}-1$ is the $(1,1)$ hook number of $\la$.
\end{prop}

\begin{pf}
In computing ${\it MurNak}(\la,\mu)$, {\it MN1inner} is invoked
precisely once for each node in the ``recursion graph". In each
one of these $r_{\la,\mu}$ invocations of {\it MN1inner}, the
length of its first parameter, $R$, is the same as the length of
the essential part of the partition sequence of $\la$, which is
$h_{1,1}(\la)+1$.

Let $t_{\textit MN1inner}(R,t)$ be the running time of {\it
MN1inner}, excluding the recursion, on input
$(R=R_1,\dots,R_s,\;t)$, $t\le k$. {\it MN1inner} performs
$\mu_t-1$ iterations in the first loop and $s-\mu_t$ iterations in
the second loop. In both loops, the time for each iteration is
$\Theta(1)$. Therefore $$ t_{\textit MN1inner}(R,t) \in \Theta(s)
.$$

It follows that each invocation of {\it MN1inner} during the
computation of {\it MurNak}($\la,\mu$) takes time
$\Theta(h_{1,1}(\la))$, with the possible exception of one
invocation with the trivial case $t>k$ which takes $\Theta(1)$.
However, this possible exception is negligible as long as
$r_{\la,\mu}>1$.

Consequently,
\begin{align*}
    t_{\textit MurNak}(\la,\mu)
        & \in t_{\textit PartSeq}(\la) + r_{\la,\mu} \Theta(h_{1,1}(\la)) \\
        & = \Theta(h_{1,1}) + r_{\la,\mu} \Theta(h_{1,1}(\la)) \\
        & = \Theta(r_{\la,\mu} h_{1,1}(\la)) \qed
\end{align*}
\end{pf}

If $\la\vdash n$ then for any $\mu\vdash n$, $\cR_{\la,\mu}
\subseteq \cR_{\la,(1^n)}=\{\alpha \st \alpha \subseteq \la \}$
and consequently $r_{\la,\mu} \le r_{\la,(1^n)}$. Define
\[r_{\la}=r_{\la,(1^n)},\] the number of partitions $\alpha$ such
that $\alpha \subseteq \la$.

\begin{lem}\label{LE:r_la}
Let $\la=(\la_1,\dots,\la_m)\vdash n$. Then \[r_\la =
\det\left({\la_i+1 \choose 1+i-j}\right)_{\substack{1\le i \le m
\\ 1 \le j \le m}}\]
\end{lem}
The lemma follows from \citet{sta1}, ch.~3, ex.~63: substituting
the empty partition for $\mu$ and $1$ for $n$, it states that
\[
\zeta^2 (\emptyset, \la) = \det\left({\la_i+1 \choose
1+i-j}\right)_{\substack{1\le i \le m
\\ 1 \le j \le m}}.
\]
By definition of $\zeta$,
\[
\zeta^2(\mu,\la) = \sum_{\mu \le \alpha \le \la} 1\quad .
\]
Noting that the partial order $\le$ defined in the exercise
coincides with containment of Young diagrams, we have
\[
\zeta^2(\emptyset,\la) = \sum_{\emptyset \le \alpha \le \la} 1 =
\sum_{\alpha \subseteq \la} 1 = \#\{ \alpha \st \alpha \subseteq
\la \} = r_\la .
\]

\subsection{\it Roich}

Let
\[
\mathcal Q_{\la,\mu} = \{ Q \st Q \text{ is a standard tableau of
shape $\alpha\subseteq\la$ and } g_\mu(Q,i)\neq 0, 1\le i \le
|\alpha| \}
\]
which is the set of standard tableaux $Q$ contained in $\la$ such
that the values \\ $g_\mu(Q,1),g_\mu(Q,2),\dots,g_\mu(Q,j)$ alone,
where $j=|\sh(Q)|$, are insufficient to determine whether
$\prod_{i=1}^{n+1}g_\mu(T,i)=0$ for all $T$ such that $T_{|j}=Q$.
(Hence, for each $Q\in\cQ_{\la,\mu}$, {\it RoiInner} is invoked to
compute $A(\la,\mu,Q)$).

For instance, in example \ref{EX:Roi},
\[
\cQ_{\la,\mu} = \left\{ \emptyset, \Dyoung{(1)}, \Dyoung{(12)},
\Dyoung{(1,2)}, \Dyoung{(13,2)}, \Dyoung{(1,2,3)},
\Dyoung{(13,2,4)}, \Dyoung{(14,2,3)} \right\}.
\]
Note that $S=\Dyoung{(12,3)} \notin \cQ_{\la,\mu}$, since $g_\mu
(S,3)=0$. Consequently, the algorithm does not compute
$A(\la,\mu,S)$.

Define \[q_{\la,\mu}=|\mathcal Q_{\la,\mu}| .\]

$\cQ_{\la,\mu}$ is just the set of non-leaf nodes in the recursion
tree of {\it RoiInner}. Since the algorithm's work on each such
node is linear in $\ell(\la)$, and the work on each leaf node is
constant, we have

\begin{prop}\label{PR:tRoi}
Let $t_{\textit{Roich}}(\la,\mu)$ be the running time of {\it
Roich} on input $(\la,\mu)$.
Then
\begin{equation}\label{EQ:tRoi}\begin{aligned}
t_{\textit{Roich}}(\la,\mu) & \in \Theta(\ell(\la) q_{\la,\mu}) .
\end{aligned}\end{equation}
\end{prop}

It is clear from \eqref{EQ:fmuQi} that if
$B(\mu)=\{1,2,\dots,n\}=B((1^n))$ then $f_\mu(Q,i) \neq 0$ for
every tableau $Q$ and $1\le i < n$. Thus for any $\mu\vdash n$,
$\mathcal Q_{\la,\mu} \subseteq \mathcal Q_{\la,(1^n)}$, and
consequently $q_{\la,\mu} \le q_{\la,(1^n)}$. Define \[q_\la =
q_{\la,(1^n)},\] the number of standard tableaux of shapes
contained in $\la$.

\begin{lem}\label{LE:q_la}
Let $\la\vdash n$. Then
\[
    d_\la \le q_\la \le n d_\la+1
\]
where $d_\la$ is the number of standard Young tableaux of shape
$\la$.
\end{lem}
\begin{pf}
Let
\[
\mathcal D_\la = \{Q \st Q\text{ is a standard tableau of shape
}\la \}.
\]
Then
\begin{align*}
\mathcal D_\la & \subseteq  \cQ_\la  \subseteq \{\emptyset \} \cup
\{ T_{|i} \st T \in \mathcal D_\la,\; 1\le i \le n \}
\end{align*}
and the lemma follows since $d_\la = |\mathcal D_\la|$.\qed
\end{pf}

\section{Comparing the Algorithms}\label{SEC:comp} \setcounter{subsection}{-1}
\subsection{Worst case analysis}

Recall that problem instances in the case of computing character
values of the symmetric group, are simply pairs $(\la,\mu)$ of
partitions of the same integer. In order to compare the
algorithms' running times, we express them as functions of the
{\it problem instance size}. A natural measure of instance size is
the weight $n$ of the partitions $\la$ and $\mu$.

In worst case analysis we consider the maximum running time of
each of the algorithms on a problem instance of size $n$, namely
\[
t_{\textit MurNak}(n) = \max_{\la,\mu \vdash n}t_{\textit
MurNak}(\la,\mu)
\]
and
\[
t_{\textit Roich}(n) = \max_{\la,\mu \vdash n}t_{\textit
Roich}(\la,\mu).
\]
By proposition \ref{PR:tMN1},
\[
t_{\textit MurNak}(n) \in \Theta(\max_{\la,\mu \vdash
n}r_{\la,\mu}h_{1,1}(\la)) = \Theta(\max_{\la\vdash n}r_\la
h_{1,1}(\la))
\]
and by proposition \ref{PR:tRoi},
\[
t_{\textit Roich}(n) \in \Theta(\max_{\la,\mu \vdash n}
\ell(\la)q_{\la,\mu}) = \Theta(\max_{\la \vdash n}
\ell(\la)q_{\la}).
\]
Hence we seek expressions for (bounds on) $\max_{\la \vdash
n}r_{\la}h_{1,1}(\la)$ and $\max_{\la \vdash n} \ell(\la)q_{\la}$
as functions of $n$.

We consider only certain families of problem instances, namely
those in which $\la$ is restricted to a given $(k,\ell)$ hook (see
below).

\subsection{$(k,\ell)$ hooks} The {\it $(k,\ell)$ hook} is the
infinite shape $\{(i,j) \st i\le k \text{ or } j \le \ell \}$. Let
$H(k,\ell;n)$ be the set of all partitions of $n$ lying inside the
$(k,\ell)$ hook, that is
$$H(k,\ell;n) = \{\la \vdash n \st \la_{k+1}\le \ell \}.$$

The two propositions in this subsections show that for partitions
in the $(k,\ell)$ hook, {\it MurNak} runs in polynomial time
$\Theta(n^{k+\ell+1})$ whereas {\it Roich}'s running time is
exponential in $n$, being in $\Omega(n^{-g} (k+\ell)^n) \cap
O(n^{-g+2} (k+\ell)^n)$ for some constant $g$.

\begin{prop}\label{PR:klhookMN}
Fix $k$ and $\ell$. Then \[\max_{\la \in H(k,\ell;n)}
h_{1,1}(\la)r_\la \in \Theta(n^{k+\ell+1})\] where
$h_{1,1}(\la)=\la_1+\la'_1 -1$.
\end{prop}

The proof requires the following lemmas. We use the notations
\begin{align*}
a \vee b &= \max \{a,b\}  \\
a \wedge b &= \min \{a,b\}.
\end{align*}

\begin{lem}\label{LE:h_11_for_hooks}
Fix $k$ and $\ell$. If $\la \in H(k,\ell;n)$ then $h_{1,1}(\la)
\in \Theta(n)$.
\end{lem}

\begin{pf}
\begin{align*}
\la & \subset \{ (i,j) \st 1 \le i \le k,\; 1 \le j \le \la_1 \}
\cup \{ (i,j) \st 1 \le j \le \ell,\; 1 \le i \le \la'_1 \}
\end{align*}
so $$n = |\la|  \le k \la_1 + \ell \la'_1 \le (k \vee
\ell)(\la_1+\la'_1) \le 2(k \vee \ell) h_{1,1}(\la) .$$ On the
other hand, $$H_{1,1}(\la)\subseteq \la$$ so
$$h_{1,1}(\la)=|H_{1,1}(\la)| \le n .$$ Therefore $$\frac{1}{2(k\vee\ell)}n \le h_{1,1}(\la) \le
n . \qed$$
\end{pf}

\begin{lem}[The $k$ strip ($\ell=0$) case]\label{LE:mRow}
Let $\la=(\la_1,\la_2,\dots,\la_k)\vdash n$. Then
\begin{equation}\label{EQ:mRow}
\frac{1}{k!}\prod_{i=1}^k (\la_i+1) \le r_\la \le \prod_{i=1}^k
(\la_i+1) .
\end{equation}
\end{lem}

\begin{pf}
Set
\[
T_\la =\{(p_1,\dots,p_k)\in \Z^k \st 0 \le p_i \le \la_i\}.
\]
Let $p=(p_1,\dots,p_k)\in T_{\la}$. Then there exists a
permutation $\sigma \in S_{k}$ such that $p_{\sigma(1)}\ge
p_{\sigma(2)}\ge \dots \ge p_{\sigma(k)}$. We claim that $\alpha =
\sigma p = (p_{\sigma (1)}, p_{\sigma (2)},\dots,p_{\sigma
(k)})\in \cR_\la$. Indeed, $\alpha$ is a partition, and for all
$i$ we have
\begin{align*}
 |\{j \st \alpha_{j} > \la_{i}\}|  =  |\{j \st p_{\sigma (j)} > \la_{i}\}|  = |\{j \st p_{j} > \la_{i}\}| <i
\end{align*}
whence $\alpha_{i} \le \la_{i}$.

It follows that for every $p\in T_\la$ there exist
$\alpha\in\cR_\la$ and $\sigma \in S_k$ such that $p=\sigma^{-1}
\alpha$. Thus $T_\la \subseteq S_k\cR_\la$, so
\begin{align*}
|T_\la| &\le |S_k|\cdot|R_\la| \\
\prod_{i=1}^k (\la_i+1) &\le k! r_\la \\
\frac{1}{k!}\prod_{i=1}^k (\la_i+1) &\le r_\la . \\
\end{align*}
The other inequality in \eqref{EQ:mRow} follows from the fact that
$\cR_\la \subseteq T_\la$.\qed
\end{pf}

\begin{pf}[Proof of proposition \ref{PR:klhookMN}]
Without loss of generality, assume $k\le \ell$. We have
\begin{equation}\label{EQ:Hkl_as_union}
H(k,\ell;n) = \bigcup_{r=1}^\ell H_i
\end{equation}
where
\[
H_i = \{ \la\in H(k,\ell;n) \st (i\wedge k,i) \in \la,\; ((i+1)
\wedge k, i+1) \notin \la \} .
\]

Let $\la \in H_i$. There is a bijection between $\cR_\la$ and
ordered triplets of partitions $(\alpha,\nu,\mu)$ such that \\
1. $\alpha$ is contained in the $(i\wedge k)\times i$ rectangle, that is $\alpha \in \cR_{(i^{i\wedge k})}$. \\
2. $\nu \subseteq (\la_1-i, \la_2-i,\dots,\la_u-i)$ where $u =
\max \{ j \st \alpha_j=i \}$. \\
3. $\mu \subseteq (\la'_1-(i\wedge k), \la'_2-(i\wedge k), \dots,
\la'_v-(i\wedge k))$ where $v = \max \{ j \st \alpha'_j= i\wedge k
\}$.\\
The following figure illustrates this bijection, showing $\la \in
H_{i}$ in white and $\rho \in \cR_\la$ shaded:

\[
\begin{picture}(0,0)%
\includegraphics{bijection.pstex}%
\end{picture}%
\setlength{\unitlength}{3947sp}%
\begingroup\makeatletter\ifx\SetFigFont\undefined
\def\x#1#2#3#4#5#6#7\relax{\def\x{#1#2#3#4#5#6}}%
\expandafter\x\fmtname xxxxxx\relax \def\y{splain}%
\ifx\x\y   
\gdef\SetFigFont#1#2#3{%
  \ifnum #1<17\tiny\else \ifnum #1<20\small\else
  \ifnum #1<24\normalsize\else \ifnum #1<29\large\else
  \ifnum #1<34\Large\else \ifnum #1<41\LARGE\else
     \huge\fi\fi\fi\fi\fi\fi
  \csname #3\endcsname}%
\else
\gdef\SetFigFont#1#2#3{\begingroup
  \count@#1\relax \ifnum 25<\count@\count@25\fi
  \def\x{\endgroup\@setsize\SetFigFont{#2pt}}%
  \expandafter\x
    \csname \romannumeral\the\count@ pt\expandafter\endcsname
    \csname @\romannumeral\the\count@ pt\endcsname
  \csname #3\endcsname}%
\fi
\fi\endgroup
\begin{picture}(5037,4218)(76,-3448)
\put(1426,-211){\makebox(0,0)[lb]{\smash{\SetFigFont{12}{14.4}{rm}$\alpha$}}}
\put(2476,-61){\makebox(0,0)[lb]{\smash{\SetFigFont{12}{14.4}{rm}$u$}}}
\put(3151,-61){\makebox(0,0)[lb]{\smash{\SetFigFont{12}{14.4}{rm}$\nu$}}}
\put(901,-961){\makebox(0,0)[lb]{\smash{\SetFigFont{12}{14.4}{rm}$v$}}}
\put(826,-1411){\makebox(0,0)[lb]{\smash{\SetFigFont{12}{14.4}{rm}$\mu'$}}}
\put(1576,694){\makebox(0,0)[lb]{\smash{\SetFigFont{12}{14.4}{rm}$i$}}}
\put( 16,-361){\makebox(0,0)[lb]{\smash{\SetFigFont{12}{14.4}{rm}$i \wedge k$}}}
\end{picture}

\]

It follows that
\[
r_\la=|\cR_\la|=\sum_{\substack{0\le u \le i\wedge k \\ 0 \le v
\le i}} a_{u,v}r_{(\la_1-i, \la_2-i,\dots,\la_u-i)}
r_{(\la'_1-(i\wedge k), \la'_2-(i\wedge k), \dots, \la'_v-(i\wedge
k))}
\]
where $$a_{u,v} = \Bigl| \{\alpha \in \cR_{(i^{i\wedge k})}\st
\max \{ j \st \alpha_j=i \}=u,\; \max \{ j \st \alpha'_j= i\wedge
k \}=v \} \Bigr|.$$ If $\alpha \in \cR_{(i^{i\wedge k})}$ and $u$,
$v$ are as above, then $u=i\wedge k \iff v=i \iff \alpha =
(i^{i\wedge k})$, thus $a_{i\wedge k,i} = 1$. Otherwise, the
partition sequence of $\alpha$ is $\dots \0 \underbrace{\1 \dots
\1}_{v \textrm{ \1s}} \0 \tilde A \1 \underbrace{\0 \dots \0}_{u
\textrm{ \0s}} \1 \dots$ where $\tilde A$ is any sequence
containing exactly $(i-v-1)$ \1s and $((i\wedge k) -u -1)$ \0s,
whence $a_{u,v} = { (i\wedge k) -u-1 + i -v-1 \choose i-v-1 }$ for
$u<i\wedge k$, $v<i$.

Thus
\begin{align*}
r_\la = |\cR_\la|  =&  \sum_{\substack{0\le u < i\wedge k \\ 0 \le
v < i}} { i\wedge k + i -u-v-2 \choose i-v-1 } r_{(\la_1-i,
\la_2-i,\dots,\la_u-i)} r_{(\la'_1-(i\wedge k), \la'_2-(i\wedge
k), \dots, \la'_v-(i\wedge k))} \\
&+ r_{(\la_1-i, \la_2-i,\dots,\la_{i\wedge k}-i)}
r_{(\la'_1-(i\wedge k), \la'_2-(i\wedge k), \dots, \la'_i-(i\wedge
k))}
\end{align*}
Since ${ i\wedge k + i -u-v-2 \choose i-v-1 }>0$ for all values of
$u$ and $v$ in the sum and does not depend on $\la$,
\begin{equation*}
\begin{split}
r_\la & \in \Theta \Big(  \sum_{\substack{1\le u \le i\wedge k \\
1 \le v \le i}}  r_{(\la_1-i, \la_2-i,\dots,\la_u-i)}
r_{(\la'_1-(i\wedge k), \la'_2-(i\wedge k), \dots, \la'_v-(i\wedge
k))}  \\  & \quad  + r_{(\la_1-i, \la_2-i,\dots,\la_{i\wedge
k}-i)} r_{(\la'_1-(i\wedge k), \la'_2-(i\wedge k), \dots,
\la'_i-(i\wedge
k))} \Big) \\
\text{(by lemma \ref{LE:mRow})} \\
& = \Theta  \Big( \sum_{\substack{1\le u \le i\wedge k \\
1 \le v \le i}}  \prod_{s=1}^u (\la_s-i+1) \prod_{t=1}^v
(\la'_t-(i\wedge k)+1) \\ & \quad + \prod_{s=1}^{i \wedge k}
(\la_s-i+1)
\prod_{t=1}^i (\la'_t-(i\wedge k)+1) \Big) \\
& = \Theta \left( \prod_{s=1}^{i \wedge k} (\la_s-i+1)
\prod_{t=1}^i (\la'_t-(i\wedge k)+1) \right)
\end{split}
\end{equation*}
Whence
\[
\max_{\la \in H_i} r_\la \in \Theta \left( \left(\frac{n-i(i\wedge
k)}{i+ (i \wedge k)} \right)^{i + (i \wedge k)} \right) =
\Theta(n^{i + (i \wedge k)})
\]
and therefore, by \eqref{EQ:Hkl_as_union},
\[
\max_{\la \in H(k,\ell;n)} r_\la = \max_{1\le i \le \ell}
\max_{\la \in H_i} r_\la  \in \Theta(n^{k+\ell})
\]
Finally, by lemma \ref{LE:h_11_for_hooks}
\[
\max_{\la \in H(k,\ell;n)} h_{1,1}(\la)r_\la \in
\Theta(n^{k+\ell+1}) \qed
\]
\end{pf}

The running time of {\it Roich} for partitions in the $(k,\ell)$
hook is determined up to a factor of order $n^2$ in the following
proposition.

\begin{prop}\label{PR:klhookRoi}
Fix $k$ and $\ell$. Then
\[
\max_{\la \in H(k,\ell;n)} \ell(\la)q_\la \in \Omega\left(
\left(\frac{1}{n}\right)^{g} (k+\ell)^n\right) \cap O\left(
\left(\frac{1}{n}\right)^{g-2} (k+\ell)^n\right)
\]
for a certain constant $g$.
\end{prop}

The proposition follows immediately from lemma \ref{LE:q_la} and
from the following theorem.

\begin{thm}[\citet{reg}, Theorem 3.3 (4)]
Assume $n$ is large and $\la \in H(k,\ell;n)$ maximizes $d_\la$.
There exist constants $c$ and $g$ such that
\[
d_\la \simeq c\left(\frac{1}{n}\right)^g (k+\ell)^n .
\]
\end{thm}

\begin{exmp}
Table \ref{TB:hook} shows the running times of the two algorithms
on $(\la,(1^{|\la|}))$ for several $\la$s in the $(1,2)$ hook.
Maximal $r_\la h_{1,1}(\la)$ and $\ell(\la) q_\la$ values for each
$n$ appear in boldface.
\end{exmp}

\subsection{General diagrams}\label{SUB:GeneralDiagrams}
Table \ref{TB:general} shows the running times of the two
algorithms for several pairs $(\la,\mu)$. The values of
$r_{\la,\mu}$ and $q_{\la,\mu}$ were obtained by running {\it
MurNak} and {\it Roich} on each pair and counting invocations of
{\it MN1Inner} and {\it RoiInner} respectively. Maximal
$r_{\la,\mu}h_{1,1}(\la)$ and $\ell(\la) q_{\la,\mu}$ values for
each $n$ appear in boldface.

\section*{Acknowledgments}
This paper is based on work conducted for my M.Sc. thesis, under
the supervision of Professor Amitai Regev. I would like to thank
him for his patient guidance, helpful advice and constant
encouragement, and specifically for reviewing and commenting on
drafts of this paper. I would also like to thank Yuval Roichman
for his comments.

\begin{table}[h]
\caption{Running times of the two algorithms on $(\la,1^{|\la|})$
when $\la$s in the (1,2) hook.}\label{TB:hook}
\begin{tabular}{|c|c||r|r||r|r|}
\cline{3-6} \multicolumn{2}{c}{} \vline &  \multicolumn{2}{c}{Murnaghan-Nakayama} \vline & \multicolumn{2}{c}{Roichman} \vline \\
\hline \parbox[c][20pt]{10pt}{$n$}  &   $\la$   & $r_{\la}$ &
$r_{\la}h_{1,1}(\la)$ & $q_{\la}$ & $\ell(\la) q_{\la}$
\\
\hline 6 &   \parbox[c][40pt]{16pt}{$\Dyng{4,1,1}$}  & 13 & {\bf 78} & 35 & 105 \\
\cline{2-6} & \parbox[c][40pt]{16pt}{$\Dyng{3,2,1}$} & 14 & 70 & 48 & {\bf 144} \\
\hline 9 & \parbox[c][40pt]{30pt}{$\Dyng{4,2,1,1,1}$} & 33 & {\bf
264} & 599 &
{\bf 2,995} \\
\hline 12 & \parbox[c][50pt]{40pt}{$\Dyng{6,2,1,1,1,1}$} & 62 &
{\bf 682} & 7,010 &
42,060 \\
\cline{2-6} & \parbox[c][50pt]{40pt}{$\Dyng{5,2,2,1,1,1}$} & 67 &
670 & 11,664 &
{\bf 69,984} \\
\hline 15 & \parbox[c][60pt]{50pt}{$\Dyng{7,2,2,1,1,1,1}$} & 116 &
{\bf 1,508} & 170,566 &
1,193,962 \\
\cline{2-6} & \parbox[c][60pt]{50pt}{$\Dyng{6,2,2,2,1,1,1}$} & 118
& 1,416 & 238,174 &
{\bf 1,667,218} \\
\hline 18 & \parbox[c][60pt]{58pt}{$\Dyng{8,2,2,2,1,1,1,1}$} & 191
& {\bf 2,865} & 4,000,428 &
32,003,424 \\
\cline{2-6} & \parbox[c][60pt]{58pt}{$\Dyng{7,2,2,2,2,1,1,1}$} &
189 & 2,646 & 5,029,991 &
{\bf 40,215,928} \\
\hline
\end{tabular}
\end{table}

\begin{table}[p]
\caption{Running times of the two algorithms for various
inputs.}\label{TB:general}
\begin{tabular}{|c|c|c||r|r||r|r|}
\cline{4-7} \multicolumn{3}{c}{} \vline &  \multicolumn{2}{c}{Murnaghan-Nakayama} \vline & \multicolumn{2}{c}{Roichman} \vline \\
\hline \parbox[c][20pt]{10pt}{$n$}  &   $\la$   &   $\mu$   &
$r_{\la,\mu}$   & $r_{\la,\mu}h_{1,1}(\la)$    &   $q_{\la,\mu}$ &
$\ell(\la) q_{\la,\mu}$
\\
\hline 6 &   \parbox[c][40pt]{16pt}{$\Dyng{3,2,1}$} & $(1^6)$ & 14 & 70 & 48 & {\bf 144} \\
\cline{3-7} & & \parbox[c][40pt]{20pt}{$\Dyng{3,2,1}$}      & 5  & 25 & 32 & 96 \\
\cline{2-7} &   \parbox[c][40pt]{16pt}{$\Dyng{3,1,1,1}$} & $(1^6)$ & 13 & {\bf 78} & 35 & 140 \\
\hline 8 & \parbox[c][40pt]{30pt}{$\Dyng{4,2,1,1}$} & $(1^8)$ & 26 & {\bf 182} & 276 & {\bf 1104} \\
\cline{3-7} & & \parbox[c][40pt]{30pt}{$\Dyng{4,2,1,1}$} &      7  & 49  & 97  & 485 \\
\hline 12 & \parbox[c][40pt]{30pt}{$\Dyng{5,3,2,1,1}$} &
$(1^{12})$ & 75 & 675 & 22,454 &
{\bf 112,270} \\
\cline{3-7} & & \parbox[c][40pt]{30pt}{$\Dyng{5,3,2,1,1}$} & 1 & 9
& 1,912 & 9,560
\\
\cline{2-7} & \parbox[c][48pt]{36pt}{$\Dyng{6,2,1,1,1,1}$} &
$(1^{12})$ & 62 & {\bf 682} & 7,010
& 42,060 \\
\cline{2-7} & \parbox[c][40pt]{30pt}{$\Dyng{4,4,2,1,1}$} &
$(1^{12})$ & 63 & 504 & 13,921
& 69,605 \\
\cline{3-7} & & \parbox[c][40pt]{30pt}{$\Dyng{4,4,2,1,1}$} & 9 &
72 & 1,384 &
6,920 \\
\hline 15 & \parbox[c][50pt]{30pt}{$\Dyng{5,4,2,2,1,1}$} & $(1^{15})$ & 139 & 1,390 &  714,201 & {\bf 4,285,206} \\
\cline{2-7} & \parbox[c][54pt]{40pt}{$\Dyng{6,3,2,1,1,1,1}$} &
$(1^{15})$ &
142 & {\bf 1,704} & 463,996 & 3,247,972 \\
\hline
\end{tabular}
\end{table}

 \end{document}